\documentclass[11pt]{article}
\usepackage{latexsym}
\usepackage{amssymb}
\usepackage{amsmath}
\usepackage{amsthm}
\usepackage{amstext}
\usepackage{amsfonts}
\usepackage[numbers]{natbib}
\usepackage{float}
\usepackage{pgf}
\usepackage{verbatim}
\usepackage{caption}
\usepackage{graphicx}
\usepackage[english]{babel}
\usepackage{color}
\usepackage{hyperref}
\usepackage{bm,mathrsfs,makeidx,graphicx}
\usepackage{subfigure}
\begin{document}
\baselineskip=.25in
\parindent=30pt

\newtheorem{theorem}{Theorem}
\newtheorem{proposition}{Proposition}
\newtheorem{lemma}{Lemma}
\newtheorem{corollary}{Corollary}
\newtheorem{definition}{Definition}
\newtheorem{example}{Example}
\newtheorem{remark}{Remark}
\newtheorem{question}{Question}

\newtheorem{tm}{Theorem}
\newcommand{\thm}{\begin{tm}}
\newcommand{\thmm}{\end{tm}}

\newtheorem{prop}{Proposition}
\newcommand{\prp}{\begin{prop}}
\newcommand{\prpp}{\end{prop}}

\newtheorem{lma}{Lemma}
\newcommand{\lm}{\begin{lma}}
\newcommand{\lmm}{\end{lma}}

\newtheorem{cro}{Corollary}
\newcommand{\cor}{\begin{cro}}
\newcommand{\corr}{\end{cro}}

\newtheorem{dfn}{Definition}
\newcommand{\df}{\begin{dfn}}
\newcommand{\dff}{\end{dfn}}

\newtheorem{exa}{Example}
\newcommand{\ex}{\begin{exa}}
\newcommand{\exx}{\end{exa}}

\newtheorem{remak}{Remark}
\newcommand{\rmk}{\begin{remak}}
\newcommand{\rmkk}{\end{remak}}

\def\bs{\vskip24pt}
\def\ms{\vskip12pt}
\def\ss{\vskip6pt}
\def\cl{\centerline}
\def\nt{\noindent}

\newcommand{\fn}{\footnote}
\newcommand{\id}{\mathop{\rm id}}

\def\qed{\quad\vrule height4pt width4pt depth0pt}

\newcommand \R { \mathbb{R} }
\newcommand \E { \mathbb{E} }
\newcommand \N { \mathbb{N} }
\newcommand \Q { \mathbb{Q} }

\def\title #1{\begin{center}
{\Large {\bf #1}}
\end{center}}
\def\author #1{\begin{center} {\large #1}
\end{center}}

\begin{titlepage}
\def\thefootnote{\fnsymbol{footnote}}
\vspace*{1.1in}
\def\date #1{\centerline {\large #1}}

\title{Conditional Exact Law of Large Numbers and Asymmetric \\
\vskip .25em
Information Economies with Aggregate Uncertainty}

\vskip .75em

\author{Lei Qiao\footnote{Department of Mathematics,
National University of Singapore, 10 Lower Kent Ridge Road,
Singapore 119076.
e-mail: a0086330@nus.edu.sg}, Yeneng Sun\footnote{Department of Economics,
National University of Singapore, 1 Arts Link, Singapore 117570.
e-mail: ynsun@nus.edu.sg}, Zhixiang Zhang\footnote{China Economics and Management Academy, Central University of Finance and Economics, Beijing 100081, China. e-mail: zhangzhixiang@cufe.edu.cn} }

\date{This version: October 3, 2014}
\vskip .75em
\vskip 1.90em
\baselineskip=.15in

\begin{abstract} A stochastic model with a continuum of economic agents often involves shocks at both macro and micro levels. This can be formalized by a continuum of random variables that are conditionally independent given the macro level shocks. Based on the framework of a Fubini extension, the results on the exact law of large numbers and its converse for a continuum of independent random variables in \cite{Sun06} are extended to the setting with conditional independence given general macro states. As an illustrative application, it is shown that any ex ante efficient allocation in an asymmetric information economy with general aggregate uncertainty has a (utility) equivalent allocation that is incentive compatible, which generalizes the corresponding results in \cite{SunYannelis07} to the case with infinitely many states.

\bs

\bs


\end{abstract}

\end{titlepage}

\setcounter{footnote}{0}

\section{Introduction}

As noted in \cite{HS08}, macroeconomic risks are the common random shocks that influence a significant portion of the population while reality suggests that these are supplemented by risks at the individual level that influence a negligible portion of the population. This can be formalized by a process with a continuum of conditionally independent random variables, given the macro level shocks. However, it is shown in Proposition 4 of \cite{HS08} that a standard joint measurability condition on the process of a continuum of random variables will imply the non-existence of individual level uncertainty. To resolve the non-compatibility of joint measurability with individual level uncertainty in our setting, we adopt the framework of a Fubini extension,\footnote{A Fubini extension is an extension of the usual product probability space that retains the Fubini property for iterated integrals. For a formal definition, see  Definition 2.2 in \cite{Sun06} and Definition \ref{df-Fubini} below.} as used in \cite{Sun06} for studying a continuum of independent random variables.

It is shown in Theorem 2.8 of \cite{Sun06} that a process measurable in a Fubini extension is essentially pairwise independent if and only if it satisfies the property of coalitional aggregate certainty in the sense that the sample distributions, based on any non-negligible collection of random variables, are essentially constant.\footnote{Theorem 7.6 of \cite{Sun98} contains the same result stated for a special type of Fubini extensions -- the Loeb product spaces.} It means that the exact law of large numbers and its converse hold in the framework of a Fubini extension. When there are only finitely many macro states, it is clear that such results are still valid in the conditional setting; see Corollary 2.11 in \cite{Sun06}. A main purpose of this paper is to show the validity of the conditional exact law of large numbers and its converse for the case with infinitely many macro states. In particular, we show in Theorem \ref{thm-CELLND} that if a process $f$ has essentially pairwise independent random variables conditioned on a countably generated $\sigma$-algebra ${\cal C}$ of events in the sample space,\fn{As noted in \cite[Example 20.1, p. 270]{pB}, a countably generated $\sigma$-algebra can be generated by a random variable taken values in $[0, 1]$.} then the sample distribution essentially equals the average regular conditional distribution of the random variables given ${\cal C}$. Furthermore, Theorem \ref{thm-CCELLND} shows the converse version of the (coalitional) conditional exact law of large numbers.

The next question is how strong the independence assumption conditioned on a countably generated $\sigma$-algebra is. For any  process $f$ measurable in a Fubini extension, it is obvious that it has event-wise measurable conditional probabilities.\footnote{For the definition of event-wise measurable conditional probabilities, see Definition 4 of \cite{HS08}.} Then Theorem 1 of \cite{HS08} shows the existence of a countably generated $\sigma$-algebra ${\cal C}$ such that the process $f$ has essentially pairwise conditionally independent random variables given ${\cal C}$. Hence, as noted formally in Lemma \ref{lm-cid} below, our assumption of conditional independence in Theorem \ref{thm-CELLND} is satisfied in general. Furthermore, basic intuition suggests that the macro shocks should reflect the uncertainty at the collective level in some appropriate sense. Theorem \ref{thm-aud} below provides some justification to this intuition by showing that the macro states can be identified with aggregate uncertainty via coalitions. In particular, the $\sigma$-algebra being conditioned on can be taken as the $\sigma$-algebra generated by the sample distributions over all the coalitions.  Theorem 2 of \cite{HS03} shows the equivalence of essential pairwise exchangeability and essential i.i.d. conditioned on some countably generated $\sigma$-algebra. Under the framework of a Fubini extension, we show in Proposition \ref{prp-ex} the equivalence of essential pairwise exchangeability and essential i.i.d. conditioned on the sample distributions.

It is well-known that there is a conflict between incentive compatibility and efficiency in a finite-agent asymmetric information economy (see, for example, \cite{GY05} p. vi, Example 0.1]). One would expect that such a conflict could be resolved in a large economy where individual agents have no monopoly power on information. This intuition has been formalized in the literature; see, for example, \cite{MP02}, \cite{MP03} and \cite{MP05} for a large but finite economy, \cite{SunYannelis07} and \cite{SunYannelis08} for an economy with a continuum of agents. In particular, \cite{SunYannelis07} considers an asymmetric information economy, where the agents have negligible information in the sense that the private signals of almost every individual can influence only a negligible group of agents and the individual agents' relevant signals are essentially pairwise independent conditioned on finitely many macro states. It is shown in \cite{SunYannelis07} that any ex ante efficient allocation in such an economy has a (utility) equivalent allocation that is incentive compatible. As noted in the above paragraph, conditional independence with possibly infinitely many macro states is generally satisfied. Thus, it is natural to ask whether the result on the consistency of incentive compatibility and ex ante efficiency in \cite{SunYannelis07} still holds when there are infinitely many macro states. By applying the conditional exact law of large numbers in our Theorem \ref{thm-CELLND}, the desired consistency result is shown in Theorem \ref{thm-XE-IC} below.

The rest of paper is organized as follows. Section \ref{sec-clln} presents the conditional exact law of large numbers and its converse, and various related results. Section \ref{sec-ic} considers the incentive compatibility problem in an asymmetric information economy with general aggregate uncertainty. All the proofs are given in the Appendix.

\section{Conditional Exact Law of Large Numbers and its Converse} \label{sec-clln}

\subsection{Basic Definitions} \label{sub-def}

Let $(I, {\cal I}, \lambda)$ be a complete atomless probability space\fn{We use the convention that all probability spaces are countably additive and complete.}, which will be the index space for the random variables in a process. Let $( \Omega, {\cal F}, P) $ be the sample probability space of the random variables in a process. In applications, the index space often represents the space of economic agents while the sample space models all the uncertainty associated with the agents. As noted in the introduction, we will use a process with a continuum of conditionally independent random variables to model both macroeconomic risks and individual level uncertainty. To resolve the non-compatibility of joint measurability with individual level uncertainty,\footnote{Proposition 4 of \cite{HS08} shows that a standard joint measurability condition on the process of a continuum of random variables will imply the non-existence of individual level uncertainty.} the framework of a Fubini extension will be used as introduced in \cite{Sun06}. Below is a formal definition of a Fubini extension as in Definition 2.2 in \cite{Sun06}. 

\df \label{df-Fubini} A probability space $(I \times \Omega, {\cal W}, Q)$ extending the usual product space $(I\times\Omega, {\cal I} \otimes {\cal F}, \lambda \otimes P)$ is said to be a Fubini extension of the usual product $(I\times\Omega, {\cal I} \otimes {\cal F}, \lambda \otimes P)$ if for any real-valued $Q$-integrable function $f$ on $(I \times \Omega, {\cal W})$,
\begin{enumerate}
\item[(1)] the two functions $f_i$ and $f_{\omega}$ are integrable, respectively, on $(\Omega, {\cal F}, P)$ for $\lambda$-almost all $i \in I$, and on $(I, {\cal I}, \lambda)$ for $P$-almost all $\omega \in \Omega$;\fn{In the sequel, we shall often use subscripts to
denote some variable of a function that is viewed as a parameter in a particular context.}
\item[(2)] $\int_{\Omega} f_idP$ and $\int_I f_{\omega}d\lambda$ are integrable, respectively, on $(I, {\cal I}, \lambda)$ and $(\Omega, {\cal F}, P)$, with $\int_{I \times \Omega}fdQ=\int_I (\int_{\Omega} f_idP) d\lambda = \int_{\Omega}(\int_I f_{\omega}d\lambda)dP$.
\end{enumerate}

\noindent To reflect the fact that the probability space $(I \times \Omega, {\cal W}, Q)$ has $(I, {\cal I}, \lambda)$ and $(\Omega, {\cal F}, P)$ as its marginal spaces, as required by the Fubini property, it will be denoted by $(I \times \Omega, {\cal I} \boxtimes {\cal F}, \lambda \boxtimes P)$.
\dff

Below we provide a formal definition of a process with conditionally uncorrelated/independent random variables based on a Fubini extension.

\df \label{df-main}
Let $(I \times \Omega, {\cal I} \boxtimes {\cal F}, \lambda \boxtimes P)$ be a Fubini extension,
${\cal C} $ be a countably generated sub-$ \sigma $-algebra of $ {\cal F}$,\footnote{Whenever necessary, we assume that a sub-$ \sigma $-algebra of $ {\cal F}$ is always strongly complete in the sense that it contains all the $P$-null sets in  ${\cal F}$.} and $X$ a complete separable metrizable topological space (i.e. a Polish space) with the Borel $\sigma$-algebra $\mathcal {B}$. Let ${\cal M} (X)$ be the space of Borel probability measures on $X$.

\begin{enumerate}
\item[(1)] Two real-valued square integrable random variables $\phi$ and $\psi$
from $( \Omega, {\cal F}, P)$ to $\mathbb{R}$
are said to be {\em conditionally uncorrelated} given ${\cal C} $ if,
the conditional expectations satisfy
\begin{equation}
     \mathbb{E}(\phi \psi |{\cal C}  )
    =  \mathbb{E}(\phi |{\cal C}  )  \mathbb{E}(\psi  |{\cal C}  ).
    \label {eq:conInd}
\end{equation}

\item[(2)] A real-valued square integrable process $f$ on $(I \times \Omega, {\cal I} \boxtimes {\cal F}, \lambda \boxtimes P)$ is said to be {\em essentially conditionally uncorrelated} given ${\cal C} $ if,
for $ \lambda $-almost all $i_1 \in I$, the random variable $ f_{ i_1 }$ is conditionally uncorrelated with $ f_{ i_2 }$
 given ${\cal C}$ for $ \lambda $-almost all $i_2 \in I$.

\item[(3)] Two random variables $\phi$ and $\psi$ from $( \Omega, {\cal F}, P)$ to  $X$ are said to be {\em conditionally independent} given ${\cal C} $ if, for any Borel sets $ B_1, B_2 \in \mathcal {B}$,
the conditional probabilities satisfy
\begin{equation}
      P(\phi ^{-1} ( B_1 )  \cap \psi ^{-1} ( B_2 ) |{\cal C}  )
    = P(\phi ^{-1} ( B_1 ) |{\cal C}  ) P(\psi ^{-1} ( B_2 ) |{\cal C}  ).
    \label {eq:conInd}
\end{equation}

\item[(4)] A process $f$ from $(I \times \Omega, {\cal I} \boxtimes {\cal F}, \lambda \boxtimes P)$ to $X$ is said to be {\em essentially pairwise conditionally independent} given ${\cal C} $ if, for $ \lambda $-almost all $i_1 \in I$,
the random variables $ f_{ i_1 }$ and $ f_{ i_2 }$ are conditionally independent given ${\cal C} $ for $ \lambda $-almost all $i_2 \in I$.\footnote{Theorem 1 of \cite{HS-ptrf} shows that essential pairwise conditional independence
is equivalent to its finite or infinite multivariate versions.}


\end{enumerate}

\dff

The following remark indicates the existence of non-trivial processes in a Fubini extension that have essentially pairwise independent random variables conditioned on any fixed countably generated $\sigma$-algebra of events.

\rmk The existence of non-trivial, independent and measurable processes in a rich Fubini extension is shown in \cite[Theorem 6.2]{Sun98} for general atomless Loeb product spaces. \cite[Proposition~5.6]{Sun06} provides another construction of a
rich Fubini extension with the unit interval $[0, 1]$ as the agent space and an extended continuum product probability space as the sample space. The main results of  \cite{SZ09} and \cite{Podcz} show respectively that the agent space can be taken as an
extended Lebesgue unit interval or a general saturated probability space.
Let $g$ be a process from a Fubini extension $(I \times \Omega, {\cal I} \boxtimes {\cal F}, \lambda \boxtimes P)$ to $\mathbb{R}$ such that the random variables $g_i, \, i \in I$ are essentially pairwise
independent. Let ${\cal C}$ be a countably generated sub-$\sigma$-algebra of ${\cal F}$; assume that it is generated by  a real valued random variable $\theta$ on $(\Omega, {\cal F}, P)$.  As in Remark 1 of \cite{HS08}, let $f$ be the process from $I \times \Omega$ to $\mathbb{R}^2$ such that $f(i, \omega) = (\theta (\omega), g(i, \omega))$ for each $(i, \omega) \in I \times \Omega$. Then $f$ is measurable in the Fubini extension $(I \times \Omega, {\cal I} \boxtimes {\cal F}, \lambda \boxtimes P)$.
By Proposition 3 in \cite{HS-ptrf}, the random variables $g_i (\cdot), \, i \in I$ are also essentially pairwise
conditionally independent given ${\cal C}$; so too are the random variables $f_i (\cdot), \, i \in I$.
\rmkk

\subsection{An Antecedent Result} \label{sub-ante}

The following exact law of large numbers is Proposition 2.5 in \cite{Sun06}.\footnote{Theorem 3.8 of \cite{Sun98} contains the same result stated for the Loeb product spaces.} It shows that the sample means of a real-valued square integrable process with essentially uncorrelated random variables in a Fubini extension are essentially constant.

\prp \label{prp-lln} Let $f$ be a real-valued square integrable process on a Fubini extension $(I \times \Omega, {\cal I} \boxtimes {\cal F}, \lambda \boxtimes P)$. If the random variables $f_i$ are essentially uncorrelated, then for $P$-almost all $\omega \in \Omega$, then the sample mean $ \mathbb{E}f_{\omega}=\int_I f_{\omega}d\lambda$ is the same as the mean of the process $ \mathbb{E}f=\int_{I \times \Omega}fd \lambda \boxtimes P$.
\prpp

Let $f'$ be a real-valued square integrable process on $(I \times \Omega, {\cal I} \boxtimes {\cal F}, \lambda \boxtimes P)$ such that the random variables $f'_i$ are essentially orthogonal with common mean zero. Then, it is obvious from the above exact law of large numbers that $\int _I f'_{\omega}(i) d\lambda=0$ for $P$-almost all $\omega \in \Omega$.

Now let $f$ be a real-valued square integrable process on $(I \times \Omega, {\cal I} \boxtimes {\cal F}, \lambda \boxtimes P)$ such that the random variables $f_i$ are essentially conditionally uncorrelated given a countably generated sub-$ \sigma $-algebra ${\cal C} $ of $ {\cal F}$. Let $g$ be the conditional expectation $ \mathbb{E}(f|{\cal I} \otimes {\cal C})$. It is easy to show that for $\lambda$-almost all $i \in I$, $g_i= \mathbb{E}(f_i|{\cal C})$ (see Lemma \ref{lm-ce} in the Appendix). Define a process $f'$ on $(I \times \Omega, {\cal I} \boxtimes {\cal F}, \lambda \boxtimes P)$ by letting $f' = f - g$. It is easy to check that the random variables $f'_i$ are essentially orthogonal with common mean zero. The following is thus an obvious corollary of Proposition \ref{prp-lln}.

\cor \label{cor-CELLN} Let $f$ be a real-valued square integrable process on $(I \times \Omega, {\cal I} \boxtimes {\cal F}, \lambda \boxtimes P)$. If $f$ is essentially conditionally uncorrelated given a countably generated sub-$ \sigma $-algebra ${\cal C} $ of $ {\cal F}$, then for $P$-almost all $\omega \in \Omega$, $\int_I f_\omega(i)d\lambda=\int_I  \mathbb{E}(f|{\cal I} \otimes {\cal C})d\lambda$.
\corr

\subsection{Conditional Exact Law of Large Numbers} \label{sub-clln}

Theorem 5.2 of \cite{Sun98} and Theorem 2.8 of \cite{Sun06} present an exact law of large numbers in terms of sample distributions. It shows that essential pairwise independence is sufficient for the sample distributions to be essentially constant. Let ${\cal C} $ be a countably generated sub-$ \sigma $-algebra of $ {\cal F}$. We shall work with essentially pairwise independent processes conditioned on ${\cal C}$. As noted in Corollary 2.11 in \cite{Sun06}, if ${\cal C}$ is generated by a finite (or countable) partition of $\Omega$, the exact law of large numbers in terms of sample distributions still holds, conditioned on each event in the partition. However, if ${\cal C}$ is not generated by a countable partition, then we can not derive the conditional exact law of large numbers by Theorem 2.8 of \cite{Sun06} directly. The purpose of this subsection is to present the conditional exact law of large numbers in the general setting.

Recall that we use a process $f$ with conditionally independent random variables given ${\cal C}$ to model both macroeconomic risks and individual level uncertainty, where ${\cal C}$ represents the macro states. In the following definition, we introduce a terminology to describe the cancelation of individual uncertainty via aggregation.

\df \label{df-nin} Let $f$ be a process from $(I \times \Omega, {\cal I} \boxtimes {\cal F}, \lambda \boxtimes P)$ to a Polish space $X$, ${\cal C} $ a countably generated sub-$ \sigma $-algebra of $ {\cal F}$, and $\mu$ a regular conditional distribution of $f$ given ${\cal I} \otimes {\cal C}$. The process $f$ is said to have no individual uncertainty in aggregation given ${\cal C}$ if $\lambda f_\omega^{-1}= \int_I\mu_{i\omega}d\lambda$ holds for $P$-almost all $\omega \in \Omega$, where $\int_I\mu_{i\omega}d\lambda$ is the random Borel probability measure $\rho$ on $X$ defined by  $\rho(B) = \int_I\mu_{i\omega}(B)d\lambda$ for any Borel set $B$ in $X$.
\dff

Let $f$, ${\cal C} $ and $\mu$ be the same as in the above definition.  The following lemma shows that $\mu$ provides the regular conditional distributions of the individual random variables in a measurable way.

\lm \label{lm-cd} Let $f$ be a process from $(I \times \Omega, {\cal I} \boxtimes {\cal F}, \lambda \boxtimes P)$ to a Polish space $X$, ${\cal C}$ a countably generated sub-$\sigma$-algebra of ${\cal F}$, and $\mu$ a regular conditional distribution of $f$ given ${\cal I} \otimes {\cal C}$.  Then, for $\lambda$-almost all $i \in I$, $\mu_i$ is a regular conditional distribution of $f_i$ given ${\cal C}$.
\lmm

We are now ready to state formally a general version of the conditional exact law of large numbers in terms of sample distributions.

\thm \label{thm-CELLND} Let $f$ be a process from $(I \times \Omega, {\cal I} \boxtimes {\cal F}, \lambda \boxtimes P)$ to a Polish space $X$. If $f$ is essentially pairwise conditionally independent given ${\cal C}$, then $f$ has no individual uncertainty in aggregation given ${\cal C}$.
\thmm

If a real-valued process $f$ is essentially pairwise conditionally independent given ${\cal C}$, then we can weaken the square integrability assumption on $f$ in Corollary \ref{cor-CELLN} to an integrability assumption on $f$, which generalizes Corollary 2.10 of \cite{Sun06} to the conditional setting.

\cor \label{cor-clln-1} Let $f$ be a real-valued integrable process on $(I \times \Omega, {\cal I} \boxtimes {\cal F}, \lambda \boxtimes P)$. If $f$ is essentially pairwise conditionally independent given ${\cal C}$, then for $P$-almost all $\omega \in \Omega$, $\int_I f_\omega(i)d\lambda=\int_I  \mathbb{E}(f|{\cal I} \otimes {\cal C})d\lambda$.
\corr

\subsection{Converse Conditional Exact Law of Large Numbers} \label{sub-cclln}


It is easy to see that for a real-valued square integrable process $f$ on a Fubini extension, the property of essentially constant sample means can not guarantee the essential uncorrelatedness for the random variables. However, part of Theorem 2.6 in \cite{Sun06} (and part of Theorem 4.6 in \cite{Sun98}) presents a converse version of exact law of large numbers in the sense that $f$ must be essentially uncorrelated if the coalitional sample means are essentially constant. The following proposition is an analog of that result in the conditional setting.

\prp \label{prp-CCELLN} Let $f$ be a real-valued square integrable process on $(I \times \Omega, {\cal I} \boxtimes {\cal F}, \lambda \boxtimes P)$. If for any $A \in \cal{I}$, $\int_A f_\omega(t)d\lambda(t)=\int_A  \mathbb{E}(f|{\cal I} \otimes {\cal C})d\lambda$ holds for $P$-almost all $\omega \in \Omega$, then $f$ is essentially conditionally uncorrelated given ${\cal C}$.
\prpp

Similarly, Theorem 2.8 in \cite{Sun06} (and part of Theorem 7.6 of \cite{Sun98}) obtains the converse exact law of large numbers in terms of sample distributions in the sense that a process on a Fubini extension must be essentially pairwise independent if the coalitional sample distributions are essentially constant. To present the converse version of the conditional exact law of large numbers in Theorem \ref{thm-CELLND}, we need the following definition first.

\df Let $f$ be a process from $(I \times \Omega, {\cal I} \boxtimes {\cal F}, \lambda \boxtimes P)$ to a Polish space $X$,
${\cal C} $ a countably generated sub-$ \sigma $-algebra of $ {\cal F}$, and $\mu$ a regular conditional distribution of $f$ given ${\cal I} \otimes {\cal C}$. The process $f$ is said to have no coalitional individual uncertainty in aggregation given ${\cal C}$ if for any $A \in {\cal I}$, and for $P$-almost all $\omega \in \Omega$, $\lambda((f^A_\omega)^{-1}(B))=\int_A\mu_{i\omega}(B)d\lambda$ holds for any Borel set $B$ in $X$, where $f^A$ is the restriction of $f$ on $A \times \Omega$.
\dff

Next we state a converse version of the conditional exact law of large numbers in terms of sample distributions.

\thm \label{thm-CCELLND} Let $f$ be a process from $(I \times \Omega, {\cal I} \boxtimes {\cal F}, \lambda \boxtimes P)$ to a Polish space $X$. If $f$ has no coalitional individual uncertainty in aggregation given ${\cal C}$, then $f$ is essentially pairwise conditionally independent given $\cal{C}$.
\thmm

\subsection{Coalitional Aggregate Uncertainty} \label{sub-au}

As indicated in the statements of Theorems \ref{thm-CELLND} and \ref{thm-CCELLND}, the assumption of conditional independence given a countably generated $\sigma$-algebra plays a key role for the validity of the conditional exact law of large numbers. An  immediate question is how restrictive when such an assumption is imposed on a general process from a Fubini extension to a Polish space. To answer this question, we need to use Theorem 1 of \cite{HS08} which provides a general characterization for a process (not necessarily measurable in a Fubini extension) to have such a structure involving conditional independence.\footnote{Theorem 1 of \cite{HS08} also indicates that the macro states could be identified via Monte Carlo simulations.}

Let $f$ be a process from $(I \times \Omega, {\cal I} \boxtimes {\cal F}, \lambda \boxtimes P)$ to a Polish space $X$, $F$ an event in ${\cal F}$, and $B$ a Borel set in $\mathcal {B}$. Then, the set $ f^{-1} (B)$ is ${\cal I} \boxtimes {\cal F}$-measurable, so is the set $D= f^{-1} (B) \cap \left(I \times F \right)$. The Fubini property implies that for $\lambda$-almost all $i \in I$, the section $D_i$ is ${\cal F}$-measurable, and the function $P(D_i) = P (F \cap f_i ^{-1} (B) )$ is $ {\cal I} $-measurable. This latter property is called event-wise measurable conditional probabilities in Definition 4 in \cite{HS08}. The following lemma is an immediate consequence of Theorem 1 of \cite{HS08}, which means that our assumption of conditional independence in Theorem \ref{thm-CELLND} is satisfied in general.

\lm \label{lm-cid}
 Let $f$ be a process from $(I \times \Omega, {\cal I} \boxtimes {\cal F}, \lambda \boxtimes P)$ to a Polish space $X$. Then, there is a countably generated $\sigma$-algebra ${\cal C}$ such that for $\lambda$-almost all $i_1 \in I$, $f_{i_1}$ and $f_{i_2}$ are conditionally independent given ${\cal C}$ for $\lambda$-almost all $i_2 \in I$.
 \lmm

Let $f$ be a process from $(I \times \Omega, {\cal I} \boxtimes {\cal F}, \lambda \boxtimes P)$ to a Polish space $X$, which has  essentially pairwise conditionally independent random variables given a countably generated $\sigma$-algebra ${\cal C}$. When $f$ is used to model both macroeconomic risks and individual level uncertainty, the individual uncertainty is reflected by the random variables $f_i$ for individual agents $i \in I$ while ${\cal C}$ represents (the information reflected by) the macro states. One can also consider the uncertainty through aggregation at the coalitional level. Let ${\underline {\cal C}}^f$ be the $\sigma$-algebra generated by the sample distributions on all the coalitions plus the $P$-null sets. That is, ${\underline {\cal C}}^f$ is generated by the ${\cal F}$-measurable mappings $\left(\{ \lambda (f^A_\omega)^{-1}(B): A\in {\cal I}, B\in{\cal B}\} \right)$ together with the $P$-null sets. Thus, ${\underline {\cal C}}^f$ reflects all the uncertainty at the aggregate level. By applying Theorem \ref{thm-CELLND} to any given coalition $A$ with $\lambda(A) >0$, it is clear that $\lambda (f^A_\omega)^{-1}(B)$ is ${\cal C}$-measurable for any $B\in{\cal B}$, and hence ${\underline {\cal C}}^f$ is always a sub-$\sigma$-algebra of ${\cal C}$.
The following theorem shows that the macro states can be identified with such coalitional aggregate uncertainty. It also implies that ${\underline {\cal C}}^f$ is the smallest sub-$\sigma$-algebra of ${\cal F}$ such that $f$ has the conditional independence structure based on such a sub-$\sigma$-algebra of ${\cal F}$.

\thm \label{thm-aud}
Let $f$ be a process from $(I \times \Omega, {\cal I} \boxtimes {\cal F}, \lambda \boxtimes P)$ to a Polish space $X$. Then the $\sigma$-algebra ${\underline {\cal C}}^f$ for the coalitional aggregate uncertainty has the following properties.

\begin{enumerate}

\item[(1)] The $\sigma$-algebra ${\underline {\cal C}}^f$ is countably generated.

\item[(2)] The process $f$ has no coalitional individual uncertainty in aggregation given ${\underline {\cal C}}^f$.


\item[(3)] The process $f$ has essentially pairwise conditionally independent random variables given ${\underline{\cal C}}^f$.
\end{enumerate}
\thmm

\subsection{Exchangeability} \label{sub-exch}

Theorem \ref{thm-aud} above shows that the macro states can be identified with the coalitional aggregate uncertainty. We shall show that for the special case of an essentially pairwise exchangeable process, the macro states can be identified with the aggregate uncertainty associated only with the grand coalition $I$. We first recall some basic definitions as in \cite{HS03} and \cite{HS08}.

\df \label{df-2} Let $f$ be a process from $(I \times \Omega, {\cal I} \boxtimes {\cal F}, \lambda \boxtimes P)$ to a Polish space $X$.
\begin{enumerate}
\item[(1)] The process $f$ is said to be {\em essentially pairwise exchangeable}
    if there exists a Borel probability measure $ \pi $
    on $ (X \times X, \mathcal {B} \otimes  \mathcal {B}) $
    such that for $ \lambda $-almost all $i_1 \in I$,  one has
   \[ P( f_{ i_1 } ^{-1} ( B_1 ) \cap f_{ i_2 } ^{-1} ( B_2 ))
   = \pi ( B_1 \times B_2 ) = \pi ( B_2 \times B_1 ) \]
for $ \lambda $-almost all $i_2 \in I$, and for all $ B_1, B_2 \in \mathcal {B}$.%

\item[(2)] Let ${\cal C}$ be a countably generated sub-$ \sigma
$-algebra of ${\cal F}$,
and let $\mu$ be a ${\cal C}$-measurable
mapping from $\Omega$ to ${\cal M} (X)$.
The process $f$ is said to be
 {\em essentially i.i.d.\ conditioned on} ${\cal C}$
if $f$ is essentially conditionally independent given ${\cal C} $,
and for $\lambda$-almost all $i \in I$,
the $ {\cal C} $-measurable mapping $ \omega \mapsto \mu _{\omega}$
is a regular conditional distribution $P( f_i^{-1} | {\cal C} )$
of the random variable $f_i$.
\end{enumerate}
\dff

Theorem 2 of \cite{HS03} and Proposition 7 of \cite{HS08} shows the equivalence of essential pairwise exchangeability and essential i.i.d. conditioned on some countably generated $\sigma$-algebra. Under the framework of a Fubini extension,the following proposition shows the equivalence of essential pairwise exchangeability and essential i.i.d. conditioned on the sample distributions.

\prp \label{prp-ex}
Let $f$ be a process from $(I \times \Omega, {\cal I} \boxtimes {\cal F}, \lambda \boxtimes P)$ to a Polish space $X$. Then $f$ is essentially pairwise exchangeable if and only if $f$ is essentially i.i.d. conditioned on the $\sigma$-algebra generated by the mapping $\lambda f_\omega^{-1}$ from $\Omega$ to ${\cal M} (X)$.

\prpp

\section{Ex ante efficiency and incentive compatibility} \label{sec-ic}

\subsection{The information structure} \label{sub-is}
\label{sub-basic}

We follow the information structure as used in \cite{SunYannelis07} and \cite{SunYannelis08} except that we use a general state space $S$ here instead of a finite state space $S$ in \cite{SunYannelis07} and \cite{SunYannelis08}. Fix an atomless probability space $(I, {\cal I}, \lambda)$ representing the space of economic agents. Let $T^0=\{q_1, q_2, \dots, q_L\}$ be the space of all the possible signals (types) for individual agents, and $\Delta^0$ the space of probability measures on $T^0$. Let $(T, {\cal T}, P^T)$ be a probability space that models the uncertainty associated with the private signal profiles for all the agents. In particular, $T$ is a space of functions from $I$ to $T^0$.\fn{In the literature, one usually assumes that different agents have possibly different sets of signals and require that the agents take all their own signals with positive probability. For
notational simplicity, we choose to work with a common set $T^0$ of signals, but allow zero probability for some of the redundant signals. There is no loss of generality in this latter approach.\label{fn-redundant} } Thus, $t \in T$, as a function from $I$ to $T^0$, represents a private signal profile for all agents in $I$. For agent $i \in I$,  $t(i)$ (also denoted by $t_i$) is the
private signal of agent $i$ while $t_{-i}$ the restriction of the signal profile $t$ to the set $I \setminus \{i\}$ of agents different from $i$; let $T_{-i}$ be  the set of all such $t_{-i}$. For simplicity, we shall assume that $(T, {\cal T})$ has a product structure so that $T$ is a product of $T_{-i}$ and $T^0$, while ${\cal T}$ is the product algebra of the power set ${\cal T}^0$ on $T^0$ with a $\sigma$-algebra ${\cal T}_{-i}$ on $T_{-i}$.  For $t \in T$ and $t_i' \in T^0$, we shall adopt the usual notation $(t_{-i}, t_i')$ to denote the signal profile whose value is $t_i'$ for agent $i$, and the same as $t$ for other agents.

Let $S$ be a complete separable metric space of true (or macro) states of nature (with its Borel $\sigma$-algebra denoted by ${\cal S}$),  which influence the utilities of all the agents, but are not known to the agents. Let $(\Omega, {\cal F}, P)$ be a probability space representing all the uncertainty on the true states as well as on the signals
for all the agents, where $(\Omega, {\cal F})$ is the product measurable space $(S \times T, {\cal S} \otimes {\cal T})$. Let
$P^S$ and $P^T$ be the marginal probability measures of $P$ respectively on $(S, {\cal S})$ and on $(T, {\cal T})$. Let $\tilde
s$ and ${\tilde t_i},\, i \in I$ be the respective projection mappings from $\Omega$ to $S$ and from $\Omega$ to $T^0$ with
${\tilde t_i}(s, t) = t_i$.\fn{${\tilde t_i}$ can also be viewed as a projection from $T$ to $T^0$.}
Denote by ${\cal C}$ the sub-$\sigma$-algebra of ${\cal F}$ generated by $\tilde s$ and the $P$-null sets. Let $(I \times \Omega, {\cal I} \boxtimes {\cal F}, \lambda \boxtimes P)$ be a Fubini extension.

For $i\in I$, let $\tau_i$ be the signal distribution of agent $i$
on the space $T^0$,\fn{For $q \in T^0$, $\tau_i (\{q\})$ is the
probability $P(\tilde t_i =q)$.} and  $P^{S\times T_{-i}}
(\cdot|t_i)$ the conditional probability measure on the product
measurable space $(S \times T_{-i}, {\cal S} \otimes {\cal T}_{-i})$
when the signal of agent $i$ is $t_i \in T^0$. If $\tau_i(\{t_i\}) >
0$, then it is clear that for $D \in {\cal S} \otimes {\cal
T}_{-i}$, $P^{S\times T_{-i}} (D|t_i) = P(D \times \{t_i\})
/\tau_i(\{t_i\})$.

Let $f$ be the private signal process from $(I \times \Omega, {\cal I} \boxtimes {\cal F}, \lambda \boxtimes P)$ to $T^0$ such that $f(i,\omega)= \tilde t_i (\omega)$, $f_\omega$ is ${\cal I}$-measurable for each $\omega \in \Omega$, and $f$ is essentially pairwise conditionally independent given ${\cal C}$.\footnote{For simplicity, here we only work with the private signal process $f$ instead of the more general idiosyncratic signal process as defined in Definition 4 of \cite{SunYannelis07}. There is no problem to use exactly the same proof to generalize our result in Theorem \ref{thm-XE-IC} to the more general setting. }

Let $\mu$ be a regular conditional distribution of $f$, given ${\cal I} \otimes {\cal C}$. As noted in Lemma \ref{lm-cd}, for $\lambda$-almost all $i \in I$, $\mu_i$ is a regular conditional distribution of $f_i$ given ${\cal C}$. Let $\bar \mu$ be the agents' average conditional signal distribution $\int_I \mu_i d\lambda(i)$. Then, the Fubini property implies that $\bar \mu$ is a ${\cal C}$-measurable mapping from $\Omega$ to $\Delta^0$. We shall impose the non-triviality
assumption on the process $f$ in the sense that the sub-$\sigma$-algebra of ${\cal F}$ generated by $\bar \mu$ and the $P$-null sets is the same as ${\cal C}$. It means that the agents' average conditional signal distribution carries as much information as the true
 states.

\subsection{The large private information economy \label{sub-model}}

We consider a large asymmetric information economy with its information structure as described in Subsection \ref{sub-is}. In this economy, agents $i \in I$ are informed with their private signals $t_i \in T^0$ but not the true states, and they can have contingent consumptions based on the signal profiles $t \in
T$ announced by all the agents. Decisions are made at the ex ante
level. The common consumption set is the positive orthant
${\mathbb R}_+^m$. In the sequel, we shall state several assumptions
on the economy.

{\it {\bf A1.} The utility function of each agent depends on her
consumption $x \in {\mathbb R}_+^m$ and the true state $s \in S$ but
not on the private signals of the agents in the economy. Thus, we
can let $u$ be a function from $I \times {\mathbb R}_+^m \times S$
to ${\mathbb R}_+$ such that for any given $i \in I$, $u(i, x, s)$
is the utility of agent $i$ at consumption bundle $x \in
{\mathbb R}_+^m$ and true state $s \in S$. For any given $x \in {\mathbb R}_+^m$, $u(i, x, s)$ is ${\cal I}\otimes {\cal S}$-measurable in $(i, s) \in I \times S$. Let $c, d \in \mathbb{R}_+$ be two constants. For any given $(i, s) \in I \times S$, $u(i, x, s)$, (also denoted by $u_{(i, s)}(x)$), is continuous and monotonic in $x \in {\mathbb R}_+^m$, and dominated by $c\|x\| + d$ for any $x \in {\mathbb R}_+^m$, where $\|\cdot\|$ is the Euclidean norm.\fn{The utility
function $u(i, \cdot, s)$ is monotonic if for any $x, y \in {\mathbb
R}_+^m$ with $x \le y$ and $x \ne y$, $u(i, x, s) < u(i, y,s)$.}

{\bf A2.} For any given $(i, s) \in I \times S$, $u(i, x, s)$, (also denoted by
$u_{(i, s)}(x)$) is concave in $x \in {\mathbb R}_+^m$.

{\bf A2$^\prime$.} For any given $(i, s) \in I \times S$, $u(i, x, s)$ is strictly concave in $x \in {\mathbb R}_+^m$.

{\bf A3.} Let $e$ be a $\lambda$-integrable function from $I$ to
${\mathbb R}_+^m$ with $e(i)$ as the initial endowment of
agent $i$.\fn{Since the true state $s \in S$ is not known to the
agents, the agents' endowments cannot depend on $s$. However, as in
\cite{MP02} and \cite{SunYannelis07}, here we also assume that the endowments
do not depend on the private signals of agents.}
}

We shall now consider an economy  where the agents are informed with
their signals but not the true state. Formally, the collection
${\cal E}^p =\{(I \times \Omega, {\cal I} \boxtimes {\cal F},
\lambda \boxtimes P), u, e,  f, \tilde s\}$ is called a Private
Information Economy.

The space of consumption plans for the economy ${\cal E}^p$ is the
space $L^1(P^T, {\mathbb R}_+^m)$ of integrable functions from $(T,
{\cal T}, P^T)$ to ${\mathbb R}_+^m$, which is infinite dimensional.
Fix an agent $i \in I$. For a consumption plan $z \in L^1(P^T,
{\mathbb R}_+^m)$, let
\begin{equation}\label{eq-U^p(z)}
U^p_i(z) = \int_{\Omega} u(i, z(t), s) dP
\end{equation}
be the ex ante expected utility of agent $i$ for the
consumption plan $z$.\fn{By assumption {\bf A1}, there are constants $c, d \in \mathbb{R}_+$ such that for any given $(i, s) \in I \times S$, $u(i, x, s) \le c\|x\| + d$ for any $x \in {\mathbb R}_+^m$, which guarantees that
$\int_\Omega u(i, z(t), s) dP$ is finite for each agent $i \in I$. \label{fn-finite-utility}}

\df \label{df-PIE}

\begin{enumerate}
\item[(1)] An allocation for the economy ${\cal E}^p$ is an integrable function $x^p$
from $(I \times T, {\cal I} \boxtimes {\cal T}, \lambda \boxtimes
P^T)$ to ${\mathbb R}_+^m$; agent $i$'s consumption plan is $x^p(i,
\cdot)$ (also denoted by $x^p_i$).

\item[(2)] An allocation $x^p$ is feasible if for $P^T$-almost all $t \in
T$, $\int_I x^p(i, t) d\lambda(i) = \int_I e(i) d\lambda(i)$.

\item[(3)] A feasible allocation $x^p$ is said to be ex ante efficient
if there does not exist a feasible allocation $y^p$ such that for
$\lambda$-almost all $i \in I$, $U^p_i(y^p_i) > U^p_i(x^p_i)$.

\item[(4)] Two allocations $x^p$ and $y^p$ are said to be utility equivalent if $U^p_i(x^p_i) = U^p_i(y^p_i)$ for $\lambda$-almost all $i \in I$.

\item[(5)] Two allocations $x^p$ and $y^p$ are said to be essentially equivalent if for $\lambda$-almost all $i \in I$, $x^p_i (t) = y^p (t)$ for $P^T$-almost all $t \in T$.

\item[(6)] For an allocation $x^p$, an agent $i \in I$, private signals $t_i, t_i' \in T^0$, let
\begin{equation} \label{eq-icaa}
U_i(x^p_i, t_i'|t_i) = \int_{S \times T_{-i}} u_i(x^p_i (t_{-i},
t_i'), s) d P^{S \times T_{-i}} (\cdot|t_i)
\end{equation}
be the interim expected utility of agent $i$ when she receives
private signal $t_i$ but mis-reports as $t_i'$. The allocation $x^p$
is said to be incentive compatible if $\lambda$-almost all $i
\in I$,
$$U_i(x^p_i, t_i|t_i) \ge U_i(x^p_i, t_i'|t_i)$$
holds for all the non-redundant signals $t_i, t_i' \in T^0$ of agent
$i$ (i.e., $\tau_i(\{t_i\}) > 0$ and $\tau_i(\{t_i'\}) > 0$).

\item[(7)] A feasible allocation $x^p$ is said to be an ex ante Walrasian
allocation (ex ante competitive equilibrium allocation) if there is
a bounded measurable price function $p$ from $(T, {\cal T})$ to
${\mathbb R}_+^m \setminus \{0\}$ such that for $\lambda$-almost all
$i \in I$, $x^p(i)$ is a maximal element in the budget set
$$
\left\{ z \in L^1(P^T, {\mathbb R}_+^m): \int_T p(t)\cdot z(t) dP^T
\le \int_T p(t)\cdot e(i) dP^T =  \left(\int_T p(t) dP^T \right)
\cdot e(i) \right\}
$$
under the expected utility function $U^p_i(\cdot)$.

\item[(8)] A coalition $A$ (i.e., a set in ${\cal I}$ with $\lambda(A)>0$) is said to ex ante
block an  allocation $x^p$ in ${\cal E}^p$ if there exists an
allocation $y^p$ such that $\int_A y^p(i, t) d\lambda(i) = \int_A
e(i) d\lambda(i)$ for $P^T$-almost all $t \in T$, and for
$\lambda$-almost all $i \in A$, $U^p_i(y^p(i)) >
U^p_i(x^p(i))$.\fn{One can also only define the allocation $y^p$ on
$A \times T$ instead of $I \times T$. However, there is no loss of
generality since one can always extend a function defined on $A
\times T$ to $I \times T$ to keep its integrability. } A feasible
allocation $x^p$ is said to be in the ex ante core of ${\cal
E}^p$, or simply an ex ante core allocation in ${\cal E}^p$, if
there is no coalition that ex ante blocks $x^p$.

\end{enumerate}
\dff

We are now ready to state two results on the general consistency between ex ante efficiency and incentive compatibility.

\thm \label{thm-XE-IC} (1) Under assumptions {\bf A1, A2, A3}, for any ex
ante efficient allocation $x^p$, there is an incentive compatible, ex ante efficient allocation $y^p$ that is utility equivalent to $x^p$.

(2) Under assumptions {\bf A1, A2$^\prime$, A3}, for any ex ante efficient allocation $x^p$, there is an incentive compatible, ex ante efficient allocation $y^p$ that is essentially equivalent to $x^p$.
\thmm

It is obvious that any ex ante core allocation is ex ante efficient.
It is also easy to check that any ex ante Walrasian allocation is ex
ante efficient. Hence the following two corollaries are clear
consequences of Theorem \ref{thm-XE-IC}.

\begin{corollary} (1) Under assumptions {\bf A1, A2, A3}, for any ex ante core allocation $x^p$, there is an incentive compatible, ex ante core allocation $y^p$ that is utility equivalent to $x^p$.

(2) Under assumptions {\bf A1, A2$^\prime$, A3}, for any ex ante core allocation $x^p$, there is an incentive compatible, ex ante core allocation $y^p$ that is essentially equivalent to $x^p$.
\end{corollary}

\begin{corollary} (1) Under assumptions {\bf A1, A2, A3}, for any ex ante Walrasian allocation $x^p$, there is an incentive compatible, ex ante Walrasian allocation $y^p$ that is utility equivalent to $x^p$.

(2) Under assumptions {\bf A1, A2$^\prime$, A3}, for any ex ante Walrasian allocation $x^p$, there is an incentive compatible, ex ante Walrasian allocation $y^p$ that is essentially equivalent to $x^p$.
 \end{corollary}

\section{Appendix} \label{sec-appendix}

\subsection{Proofs of the results in Subsection \ref{sub-clln}} \label{sub-proof-clln}

In Corollary \ref{cor-CELLN}, we work with $ \mathbb{E}(f|{{\cal I} \otimes {\cal C}})$ to state a version of the  conditional exact law of large numbers. The following lemma shows that $ \mathbb{E}(f|{{\cal I} \otimes {\cal C}})$ provides the conditional expectations of individual random variables in a measurable way.

\lm \label{lm-ce}
Let $f$ be a real-valued integrable process on $(I \times \Omega, {\cal I} \boxtimes {\cal F}, \lambda \boxtimes P)$,  ${\cal C}$ a countably generated sub-$\sigma$-algebra of ${\cal F}$, and  $g =  \mathbb{E}(f|{{\cal I} \otimes {\cal C}})$.  Then, for $\lambda$-almost all $i \in I$, the conditional expectation $ \mathbb{E}(f_i|{\cal C}) = g_i.$
\lmm

\nt{\bf Proof:} \ Fix any $C \in {\cal C}$. By the definition of conditional expectation, the identity
$$\int \!\!\int _{A\times C} f \; d\lambda \boxtimes P = \int \!\! \int_{A \times C} gd\lambda \boxtimes P$$
holds for any $A \in {\cal I}$.  Hence, the Fubini property implies that
$$\int_A \left[ \int_C (f_i(\omega) - g_i(\omega))dP\right] d\lambda = 0.$$
Thus, for any $C \in {\cal C}$, we have
$\int_C (f_i(\omega) - g_i(\omega))dP = 0$
for $\lambda$-almost all $i \in I$.

Take a countable sub-collection $\{C_n\}^\infty_{n=1}$ of ${\cal C}$ such that it generates ${\cal C}$ and is closed under finite intersections.  The previous paragraph shows that for each $n \ge 1$,
$\int_{C_n}(f_i - g_i)dP = 0$
holds for all $i \in I_n$ with $\lambda(I_n) = 1$.  Let $I_0 = {\mathop{\bigcap}^\infty_{n=1}} I_n$.  Then $\lambda(I_0) = 1$, and for each $i \in I_0$,$\int_{C_n}(f_i-g_i)dP = 0$
for each $n \ge 1$. Since the $C_n$'s form a $\pi$-system, the well-known $\pi-\lambda$ theorem\footnote{See, for example, Theorem 3.2 in p42 of \cite{pB}.}  implies that
 for each $i \in I_0$, $\int_C (f_i-g_i)dP = 0$ for all $C \in {\cal C}$.

By the classical Fubini theorem, there is $A_0 \subseteq I$ with $\lambda(A_0) = 1$ such that $g_i$ is ${\cal C}$-measurable for each $i \in A_0$.  Let $A_1 = I_0 \cap A_0$.  Then $\lambda(A_1) = 1$, and for each $i \in A_1$, $g_i$ is ${\cal C}$-measurable, and $\int_C f_i \; dP = \int_C g_i \; dP$ for all $C \in {\cal C}$.  Hence, for each $i \in A_1$, $ \mathbb{E}(f_i|{\cal C}) = g_i$. \qed

\ms

\nt{\bf Proof of Lemma \ref{lm-cd}:} By the Fubini property, we know that for $\lambda$-almost all $i \in I$, $\mu_i$ is ${\cal C}$-measurable.
Let $\{O_n\}_{n \ge 1}$ be a countable base of open sets in $X$ which is closed under finite intersections.  Then, by Corollary 2 on \cite[p.227]{CT}, we have for each $n \ge 1$,
$$
 \mathbb{E}( {\bf 1}_{O_n}(f)|{{\cal I} \otimes {\cal C}} )(i,\omega)
= \int_X {\bf 1}_{O_n}(x)d \mu(i, \omega)(x) =\mu(i,\omega)(O_n).
$$
By Lemma \ref{lm-ce}, there exists a set $I_0 \in {\cal I}$ with $\lambda(I_0) =1$ such that for any $i \in I_0$, $\mu_i$ is ${\cal C}$-measurable, and $P\left( f_i^{-1}(O_n)|{\cal C} \right)=  \mathbb{E}({\bf 1}_{O_n}(f_i)|{\cal C}) = \mu_i(\omega)(O_n)$ for all $n \ge 1$. This implies that for any $i \in I_0$, $\mu_i$ is a regular conditional distribution of the random variable $f_i$ given ${\cal C}$. \qed

\ms
\nt{\bf Proof of Theorem \ref{thm-CELLND}:}  Let $\mu$ be a regular conditional distribution of $f$ given ${\cal I} \otimes {\cal C}$. Let $\{O_n\}_{n \ge 1}$ be a countable base of open sets in $X$ which is closed under finite intersections. Then, for each $n \ge 1$, we have
\begin{equation} \label{eq-rcd}
 \mathbb{E}( {\bf 1}_{O_n}(f)|{{\cal I} \otimes {\cal C}} )(i,\omega)
= \int_X {\bf 1}_{O_n}(x)d \mu(i, \omega)(x) = \mu(i,\omega)(O_n).
\end{equation}
By the conditional independence assumption, we know that for each $n \ge 1$, for $\lambda$-almost all $i_1 \in I$,
$$P(f^{-1}_{i_1}(O_n) \cap f^{-1}_{i_2}(O_n)|{\cal C}) = P(f^{-1}_{i_1}(O_n)|{\cal C}) \cdot P(f^{-1}_{i_2}(O_n)|{\cal C})$$
holds for $\lambda$-almost all $i_2 \in I$.  Thus, for any given $n \ge 1$, the random variables in the process ${\bf 1}_{O_n}(f)$ are essentially uncorrelated conditioned on ${\cal C}$. It then follows from Corollary \ref{cor-CELLN} and Equation (\ref{eq-rcd}) that for $P$-almost all $\omega \in \Omega$,
$$
\lambda(f^{-1}_\omega(O_n))= \int_I {\bf 1}_{O_n}(f_\omega(i))d\lambda(i) = \int_I \mu(i,\omega)(O_n) d\lambda(i)
$$
holds for all $n \ge 1$. Therefore, $\lambda f^{-1}_\omega = \int_I \mu_\omega(i)d\lambda$ for $P$-almost all $\omega \in \Omega$. That is, $f$ has no individual uncertainty in aggregation given ${\cal C}$. \qed

\subsection{Proofs of the results in Subsection \ref{sub-cclln}} \label{sub-proof-cclln}

For a given real-valued integrable process $f$ on a Fubini extension $(I \times \Omega, {\cal I} \boxtimes {\cal F}, \lambda \boxtimes P)$, and a countably generated $\sigma$-algebra ${\cal C}$ of events, one can take the conditional expectations of $f$ given ${{\cal I} \otimes {\cal C}}$ and ${{\cal I} \otimes {\cal F}}$ respectively. The following lemma characterizes when the two conditional expectations are equal.

\lm \label{lm-cei} Let $f$ be a real-valued integrable process on $(I \times \Omega, {\cal I} \boxtimes {\cal F}, \lambda \boxtimes P)$,  ${\cal C}$ a countably generated sub-$\sigma$-algebra of ${\cal F}$, $g =  \mathbb{E}(f|{{\cal I} \otimes {\cal C}})$ and $h^* =  \mathbb{E}(f|{{\cal I} \otimes {\cal F}})$. Then the following are equivalent.
\begin{enumerate}
\item [(1)] For any fixed $A \in {\cal I}$, $\int_A f_\omega(i)d\lambda=\int_A  \mathbb{E}(f|{\cal I} \otimes {\cal C})d\lambda$ for $P$-almost all $\omega \in \Omega$.

\item [(2)] $h^* = g$.

\item [(3)] $h^*_i$ is $\cal C$-measurable for $\lambda$-almost all $i \in I$.
\end{enumerate}

\lmm

\nt{\bf Proof:} First consider $(2) \Longrightarrow (1)$. Assume (2), i.e., $h^* = h$, and fix any $A \in {\cal I}$.  Then, for any $F \in {\cal F}$, we have
$$\int\!\!\int_{A\times F}f \; d\lambda \boxtimes P = \int\!\!\int_{A\times F} \mathbb{E}(f|{\cal I} \otimes {\cal F})d\lambda \boxtimes P = \int\!\!\int_{A\times F}g \; d\lambda \boxtimes P.$$
Hence by Fubuni property, for any $F \in {\cal F}$,
$$\int_F \int_A f_\omega d\lambda \; dP = \int_F \int_A g_\omega \; d\lambda \; dP,$$
which implies that for $P$-almost all $\omega \in \Omega$, $\int_A f_\omega d\lambda = \int_A g_\omega d\lambda$.

\ms
Next, assume that (1) holds, i.e., for any fixed $A \in {\cal I}$, $\int_A f_\omega d\lambda = \int_A g_\omega d\lambda$ holds for $P$-almost all $\omega \in \Omega$. For any $F \in {\cal F}$, $\int_F \int_A f_\omega d\lambda \; dP = \int_F \int_A h_\omega d\lambda \; dP.$ Thus, $\int\!\!\int_{A\times F}  \mathbb{E}(f|{\cal I} \otimes {\cal F})d\lambda \boxtimes P = \int\!\!\int_{A\times F} g \; d\lambda \otimes P$.  By the $\pi-\lambda$ theorem, $h^* = g$, i.e., (2) holds.

$(2) \Longrightarrow (3)$ follows from the classical Fubini theorem.

It remains to show $(3) \Longrightarrow (2)$. Assume (3). It is clear that $ \mathbb{E}(h^* |{{\cal I} \otimes {\cal C}}) = g$. By Lemma \ref{lm-ce}, we have for $\lambda$-almost all $i \in I$,
$ \mathbb{E}(h^*_i |{\cal C}) = h_i$. Since $h^*_i$ is $\cal C$-measurable for $\lambda$-almost all $i \in I$, we have $h^*_i = g_i$ for $\lambda$-almost all $i \in I$. Applying the Fubini property again, we have for any $A \in {\cal I}$ and $F \in {\cal F}$,
$$\int\!\!\int_{A\times F}h^* d\lambda \otimes P=\int_A\int_F h^*_i dP d\lambda  = \int_A\int_F g_i dP d\lambda = \int\!\!\int_{A\times F} g d \lambda \otimes P,$$
which implies (2) by the $\pi-\lambda$ theorem. \qed








\ms

\nt{\bf Proof of Proposition \ref{prp-CCELLN}:} \ Assume that for any fixed $A \in {\cal I}$, $\int_A f_\omega(i)d\lambda=\int_A  \mathbb{E}(f|{\cal I} \otimes {\cal C})d\lambda$ for $P$-almost all $\omega \in \Omega$.
Let $g= \mathbb{E}( f | {{\cal I}\otimes {\cal C}} )$ and $h^* =  \mathbb{E}(f|{{\cal I} \otimes {\cal F}})$. Lemma \ref{lm-cei} implies that  $g=h^*$. Since $ \mathbb{E}(f -h^* | {\cal I} \otimes {\cal F})=0$, we have
$ \mathbb{E}(f -g | {\cal I} \otimes {\cal F})=0$.

Let $\alpha$ be any real-valued square integrable random variable on $(\Omega, {\cal F}, P)$. We have
$  \mathbb{E}( (f -g ) \alpha | {\cal I} \otimes {\cal F})=0.$  Taking the conditional expectation further on ${\cal I}\otimes {\cal C}$, we obtain $ \mathbb{E}( (f -g ) \alpha | {\cal I}\otimes {\cal C})=0.$ It follows from Lemma \ref{lm-ce} that $ \mathbb{E}( (f_i-g_i)\alpha | {{\cal C}}) =0$ for $\lambda$-almost all $i \in I$. Hence, by applying Lemma \ref{lm-ce} again, we obtain that $\alpha$ and $f_i$ are conditionally uncorrelated given ${\cal C}$ for $\lambda$-almost all $i \in I$. To show $f$ is essentially conditionally uncorrelated given ${\cal C}$, we just need to fix $\alpha=f_{i'}$ for $\lambda$-almost all $i' \in I$. \qed


\lm \label{lm-orth} Let $f$ and $f'$ be any real-valued square integrable processes on $(I \times \Omega, {\cal I} \boxtimes {\cal F}, \lambda \boxtimes P)$,  and ${\cal C}$ a countably generated sub-$\sigma$-algebra of ${\cal F}$. Suppose  $ \mathbb{E}(f|{\cal I} \otimes {\cal F}) = 0$. Then we have for $\lambda$-almost all $i_1 \in I$, $ \mathbb{E}(f'_{i_1} f_{i_2}|{\cal C}) = 0$ for $\lambda$-almost all $i_2 \in I$.
\lmm

\nt{\bf Proof:} \ For $\lambda$-almost all $i_1 \in I$, $f'_{i_1}$ is square integrable on $(\Omega, {\cal F}, P)$; fix such an $i_1 \in I$. By taking $\alpha=f'_{i_1}$ and $g=0$ as in the proof of Proposition \ref{prp-CCELLN}, we obtain that $ \mathbb{E}(f'_{i_1} f_{i_2}|{\cal C}) = 0$ for $\lambda$-almost all $i_2 \in I$. \qed

\ms

\nt{\bf Proof of Theorem \ref{thm-CCELLND}:} Let $\mu$ be a regular conditional distribution of $f$ given ${\cal I} \otimes {\cal C}$ and $\{O_n\}_{n \ge 1}$ be a countable base of open sets in $X$ which is closed under finite intersections. Fix $n \ge 1$, for any $A \in {\cal I}$, the assumption of no coalitional individual uncertainty in aggregation implies that
for $P$-almost all $\omega \in \Omega$, $\lambda((f_\omega^A)^{-1}(O_n))=\int_A\mu_{i\omega}(O_n)d\lambda$, and thus it follows from Lemma \ref{lm-cd} that
\begin{equation} \label{eq-cllncd}
\int_A {\bf 1}_{f^{-1}_{\omega}(O_n)}{(i)} d\lambda=\int_A  \mathbb{E}({\bf 1}_{f^{-1}_{\omega}(O_n)}(i)|{\cal I} \otimes {\cal C})d\lambda.
\end{equation}
Let $f^n(i,\omega)={\bf 1}_{f^{-1}(O_n)}(i, \omega)$ and $g^n= \mathbb{E}({\bf 1}_{f^{-1}(O_n)}|{\cal I} \otimes {\cal C}))$. By Equation (\ref{eq-cllncd}) and Lemma \ref{lm-cei}, we have $ \mathbb{E}(f^n |{\cal I} \otimes {\cal F})=g^n$, which means $ \mathbb{E}(f^n -g^n |{\cal I} \otimes {\cal F})=0$.
For any given $m \ge 1$, Lemma \ref{lm-orth} implies that for $\lambda$-almost all $i_1 \in I$, ${\mathbb E} \left( {\bf 1}_{f^{-1}_{i_1}(O_m)}{\omega)} \cdot [f^n_{i_2}-g^n_{i_2}] |{\cal C} \right) = 0$ for $\lambda$-almost all $i_2 \in I$. This means that for $\lambda$-almost all $i_1 \in I$,
$$P(f^{-1}_{i_1}(O_m) \cap f^{-1}_{i_2}(O_n)|{\cal C}) = P(f^{-1}_{i_1}(O_m)|{\cal C})P(f^{-1}_{i_2}(O_n)|{\cal C})$$
holds for $\lambda$-almost all $i_2 \in I$.  By grouping countably many null sets together as in the proof of Theorem 2.8 of \cite{Sun06}, we can obtain that for $\lambda$-almost all $i_1 \in I$,
$$P(f^{-1}_{i_1}(O_m) \cap f^{-1}_{i_2}(O_n)| {\cal C}) =
P(f^{-1}_{i_1}(O_m)| {\cal C})P(f^{-1}_{i_2}(O_n)| {\cal C})$$
holds $\lambda$-almost all $i_2 \in I$, and for all $m,n \ge 1$.
By Theorem 1 of \cite[p.230]{CT}, $f_{i_1}$ and $f_{i_2}$ are essentially pairwise conditionally independent given ${\cal C}$. \qed

\subsection{Proofs of the results in Subsections \ref{sub-au} and \ref{sub-exch}}

To prove Theorem \ref{thm-aud}, we shall first work with a real-valued square integrable process $f$ on $(I \times \Omega, {\cal I} \boxtimes {\cal F}, \lambda \boxtimes P)$, and consider an analog of Theorem \ref{thm-aud} in the setting of coalitional sample means. Let ${\cal C}_0^f$ be the $\sigma$-algebra generated by the ${\cal F}$-measurable mappings $\{ \int_A f_\omega(i) d\lambda: A\in {\cal I} \}$ together with the $P$-null sets. The following proposition presents some properties about ${\cal C}_0^f$.

\prp \label{prp-aue}
Let $f$ be a real-valued square integrable process on $(I \times \Omega, {\cal I} \boxtimes {\cal F}, \lambda \boxtimes P)$.   Then we have following properties.

\begin{enumerate}

\item[(1)] The $\sigma$-algebra ${\cal C}_0^f$ is countably generated.

\item[(2)] For any $A \in {\cal I}$, $\int_A f_\omega(i)d\lambda=\int_A  \mathbb{E}(f|{\cal I} \otimes {\cal C}_0^f)d\lambda$ for $P$-almost all $\omega \in \Omega$.


\item[(3)] The process $f$ is essentially conditionally uncorrelated given ${\cal C}_0^f$.
\end{enumerate}
\prpp

\nt{\bf Proof:} By Lemma \ref{lm-cid}, there exists a countably generated $\sigma$-algebra ${\cal C}$ such that $f$ is essentially pairwise conditionally independent given ${\cal C}$. Hence, $f$ is also essentially conditionally uncorrelated given ${\cal C}$.
Without loss of generality, we work with the strong completion of such a countably generated $\sigma$-algebra. Fix any $A \in {\cal I}$. Then Corollary \ref{cor-CELLN} implies that $\int_A f_\omega d\lambda$ is ${\cal C}$-measurable. By the arbitrary choices of $A$, we obtain that $ {\cal C}_0^f \subseteq {\cal C}$. This means that (1) holds.

Fix any $A \in {\cal I}$. The Fubini property implies that for any $C \in {\cal C}^f_0$,
$$
\int_C \int_A f_\omega d\lambda dP=\int\!\!\int_{A\times C} f d\lambda\boxtimes P
= \int\!\!\int_{A\times C}  \mathbb{E}(f|{\cal I}\otimes {\cal C}_0^f)d\lambda\boxtimes P
=\int_C\int_{A}  \mathbb{E}(f|{\cal I}\otimes  {\cal C}^f_0)d\lambda d P.
$$
Since both $\int_A f_\omega d\lambda$ and $\int_{A}  \mathbb{E}(f|{\cal I}\otimes  {\cal C}^f_0)d\lambda$ are ${\cal C}_0^f$-measurable, the above equation implies that $\int_A f_\omega d\lambda = \int_{A}  \mathbb{E}(f|{\cal I}\otimes  {\cal C}^f_0)d\lambda$,
 which means that (ii) holds.

(iii) follows from Proposition \ref{prp-CCELLN}. \qed

\ms

\nt{\bf Proof of Theorem \ref{thm-aud}:} Let $f$ be a process from $(I \times \Omega, {\cal I} \boxtimes {\cal F}, \lambda \boxtimes P)$ to a Polish space $X$. As above Lemma \ref{lm-cid} implies the existence of a countably generated $\sigma$-algebra ${\cal C}$ such that for $\lambda$-almost all $i_1 \in I$, $f_{i_1}$ and $f_{i_2}$ are conditionally independent given ${\cal C}$ for $\lambda$-almost all $i_2 \in I$. Without loss of generality, we work with the strong completion of such a countably generated $\sigma$-algebra. Fix any $A \in {\cal I}$. Then Theorem \ref{thm-CELLND} implies that for any $B \in {\cal B}$, $\lambda \left( \left( f^A_\omega \right)^{-1}(B)\right)$ is ${\cal C}$-measurable. By the arbitrary choices of $A$ and $B$, we obtain that ${\underline {\cal C}}^f \subseteq {\cal C}$. This means that (1) holds.

Let $\mu$ be  a regular conditional distribution of $f$ given ${\cal I} \otimes {\underline {\cal C}}^f$. Fix any $A \in {\cal I}$. For any $B \in {\cal B}$, $\lambda\left( \left( f_\omega\right)^{-1}(B)\right)$ is ${\underline {\cal C}}^f$-measurable. By the same idea as in the proof of Proposition \ref{prp-aue} (2), we can obtain that $\lambda\left( \left( f_\omega\right)^{-1}(B)\right) = \int_A \mu_{i \omega} (B) d\lambda$. Thus (2) holds.

(3) follows from Theorem \ref{thm-CCELLND} and parts (1) and (2) of this theorem. \qed

\ms

\nt{\bf Proof of Proposition \ref{prp-ex}:} Suppose that $f$ is essentially pairwise exchangeable. Proposition 7 of \cite{HS08} shows the existence of a ${\cal F}$-measurable mapping $\mu$ from $\Omega$ to ${\cal M} (X)$ such that $f$ is essential i.i.d. conditioned on the countably generated $\sigma$-algebra ${\cal C}$ generated by $\mu$, and $\mu$ is a regular conditional distribution of $f_i$ given ${\cal C}$ for $\lambda$-almost $i \in I$.

Let $\mu'$ be a regular conditional distribution of $f$ given ${\cal I} \otimes {\cal C}$. By Lemma \ref{lm-cd}, $\mu'_i$ is a regular conditional distribution of $f_i$ given ${\cal C}$ for $\lambda$-almost $i \in I$. Hence, for $\lambda$-almost $i \in I$, $\mu'_i = \mu$. Theorem \ref{thm-CELLND} implies that or $P$-almost all $\omega \in \Omega$, $\lambda f^{-1}_\omega = \int_I \mu'_{i\omega } d\lambda = \int_I \mu_{\omega } d\lambda = \mu_\omega$. Therefore, ${\cal C}$ is generated by the sample distribution $\lambda f^{-1}_\omega$. \qed


\subsection{Proofs of the results in Section \ref{sec-ic}} \label{sec-proof-ic}

\nt{\bf Proof of Theorem \ref{thm-XE-IC}:} (1) Let $x^p$ be any ex ante efficient allocation. For any $\omega \in \Omega$, denote the realized signal distribution $\lambda f_{\omega}^{-1}$ by $\gamma(\omega)$;
thus $\gamma$ is a ${\cal F}$ measurable mapping from $\Omega$ to $\Delta^0$. Note that the definition of the signal process $f$ does not depend on $s \in S$, and hence $\lambda f_{\omega}^{-1}$ can also be viewed as a ${\cal T}$-measurable function from $T$ to $\Delta^0$. By Theorem \ref{thm-CELLND}, we have  $\gamma(\omega)= \lambda f_{\omega}^{-1} = \bar \mu(\omega)$ for
$P$-almost all $\omega \in \Omega$. The non-triviality assumption implies that the sub-$\sigma$-algebra of ${\cal F}$ generated by $\gamma$ and the $P$-null sets is the same as ${\cal C}$, which implies that there is a Borel measurable function  $\beta$ from $\Delta^0$ to $S$ such that for $P$-almost all $\omega \in \Omega$, $\tilde s (\omega) = \beta (\gamma (\omega))$.

For any $i \in I$ and $\omega = (s, t) \in \Omega$, let $\bar{x}^p(i,\omega)=x^p(i,t)$. We shall use $\rho$ to denote the regular conditional distribution $\lambda \boxtimes P\left((\bar x^p)^{-1}|{\cal I} \otimes {\cal C}\right)$ of $\bar x^p$ given ${\cal I} \otimes {\cal C}$.
Then there exists a ${\cal I} \times {\cal S}$-measurable function $\nu$ from $I \times S$ to ${\cal M}(\mathbb{R}^+_M)$ such that $\rho(i,\omega)=\nu(i, \tilde{s}(\omega))$ for $\lambda \otimes P$-almost all $(i, \omega) \in I \times \Omega$. Let $\alpha$ be the ${\cal I} \times {\cal S}$-measurable mapping from $I \times S$ to ${\mathbb R}_+^m$ defined by $\alpha (i, s) = \int_{\mathbb{R}^m_+}x \; d\nu(i, s)$. Thus, we have for $\lambda \otimes P$-almost all $(i, \omega) \in I \times \Omega$,
$$
\mathbb{E}(\bar{x}^p|{\cal I} \otimes {\cal C})=\int_{\mathbb{R}^m_+}x \; d\rho(i, \omega)
= \int_{\mathbb{R}^+_M}x \; d\nu(i, {\tilde s}(\omega)) = \alpha(i, {\tilde s}(\omega)) = \alpha(i, \beta(\gamma(\omega))).
$$

Since $\gamma$ can be viewed as a ${\cal T}$-measurable function from $T$ to $\Delta^0$, $\alpha (i, \beta(\gamma(t)))$ defines  a ${\cal I} \otimes {\cal T}$-measurable function from $I \times T$ to ${\mathbb R}_+^m$.
For each $(i, t) \in I \times T$, let $y^p(i, t) = \alpha (i, \beta(\gamma(t))).$ It is clear that $y^p$ is ${\cal I} \otimes {\cal T}$-measurable and hence ${\cal I} \boxtimes {\cal T}$-measurable.
For $\lambda \otimes P$-almost all $(i, \omega) \in I \times \Omega$, $\bar{y}^p(i,\omega) = y^p(i, ({\tilde t}(\omega))) = \alpha(i, \beta(\gamma(\omega))) = \mathbb{E}(\bar{x}^p|{\cal I} \otimes {\cal C})$.
By Lemma \ref{lm-ce}, we obtain that for $\lambda$-almost all $i \in I$, $\bar{y}^p_i (\omega) = {\mathbb E} \left( \bar{x}_i^p | {\cal C} \right) = \alpha_i (\beta(\gamma(\omega)))$ for
$P$-almost all $\omega \in \Omega$
.
For any $C \in {\cal C}$, the Fubini property implies that
\begin{eqnarray}
\int_C\int_I y^p(i,\tilde{t}(\omega))d\lambda d P &=&
\int_{I\times C} y^p(i,\tilde{t}(\omega))d\lambda \otimes P=
\int_{I \times C} {\mathbb E} \left( \bar{x}^p | {\cal I} \otimes {\cal C} \right)d\lambda \otimes P \nonumber \\
&=& \int_{I \times C}  \bar{x}^p d\lambda \otimes P = \int_C\int_I \bar{x}^p d\lambda d P = \int_C \int_{I} e(i)d\lambda dP.
\end{eqnarray}
By the arbitrary choice of $C \in {\cal C}$, we obtain that $\int_I y^p(i,\tilde{t}(\omega))d\lambda = \int_{I} e(i)d\lambda$ for $P$-almost all $\omega \in \Omega$. Hence $\int_I y^p(i,t)d\lambda = \int_{I} e(i)d\lambda$ for $P^T$-almost all $t \in T$,
which means that $y^p$ is a feasible allocation in ${\cal E}^p$.

Next we show that the allocation $y^p$ is incentive compatible. Since no single agent could change the realized signal distribution, it is obvious that for any $i \in I$, $\gamma(t_{-i},t_i)=\gamma(t_{-i},t_i')$ for all $t_{-i}\in T_{-i}$ and all $ t_i, t_i'\in T^0$. Hence, for $\lambda$-almost all $i \in I$, for any signal $t_i$ with positive probability, we obtain that
\begin{eqnarray*}
U_i(y^p_i, t_i'|t_i) &=& \int_{S \times T_{-i}} u_i(y^p_i (t_{-i},
t_i'), s) d P^{S \times T_{-i}} (\cdot|t_i)\\
&=& \int_{S \times T_{-i}} u_i(\alpha_i (\beta (\gamma(t_{-i},
t_i'))), s) d P^{S \times T_{-i}} (\cdot|t_i)\\
&=& \int_{S \times T_{-i}} u_i(\alpha_i (\beta(\gamma(t_{-i},
t_i))), s) d P^{S \times T_{-i}} (\cdot|t_i)\\
&=&U_i(y^p_i, t_i|t_i).
\end{eqnarray*}
Therefore, the allocation $y^p$ is incentive compatible.

For $i \in I$, $\omega = (s, t) \in \Omega$, let $v(i, \omega)=u(i,x^p_i(t),s)$. Then, we have
$$\mathbb{E}(v|{\cal I} \times {\cal C})=\int_{\mathbb{R}^m_+}u(i,x,\tilde{s}(\omega))d\rho(i, \omega).$$
It follows from Lemma \ref{lm-ce} that for $\lambda$-almost all $i \in I$, $\mathbb{E}(v_i|{\cal C}) = \mathbb{E}(v|{\cal I} \times {\cal C})$. Since $u(i, \cdot, s)$ is concave, Jensen's inequality (see Lemma 3.5 in \cite[p.49]{Kalle}) implies that for $\lambda$-almost all $i \in I$
\begin{eqnarray} \label{eq-Jensen}
U_i(x^p_i) &=& \int_{\Omega} v(i,\omega) d P = \int_{\Omega} \mathbb{E}(v_i|{\cal C}) d P
= \int_{\Omega} \mathbb{E} \left(v(i,\omega)|{\cal I} \otimes {\cal C} \right) d P\nonumber \\
& = &\int_{\Omega}\int_{\mathbb{R}^m_+}u(i,x,\tilde{s}(\omega)) d \rho(i, \omega)dP\nonumber\\
& \leq &\int_{\Omega} u(i,\int_{\mathbb{R}^m_+}x \; d\rho(i, \omega), {\tilde s}(\omega)) dP\nonumber\\
& = &  \int_{\Omega}  u \left(i, \mathbb{E}(\bar{x}^p|{\cal I} \otimes {\cal C}), {\tilde s}(\omega) \right)  d P\nonumber\\
&=& \int_{\Omega} u_i(\bar y^p_i, {\tilde s}(\omega)) d P = U_i (y^p_i).
\end{eqnarray}

Let $d_i=U_i(y^p_i)-U_i(x^p_i)$ for each $i\in I$, and $A=\{i\in I: d_i>0\}$. If $\lambda(A)>0$, then
\begin{eqnarray} \label{eq-ef}
&&\int_A \int_T u(i,y^p_i(t),\beta(\gamma(t)))dP^Td\lambda = \int_A U_i (y^p_i) d\lambda > \int_A U_i (x^p_i) d\lambda \nonumber \\
&&= \int_A \int_T u(i,x^p_i(t),\beta(\gamma(t))) dP^Td\lambda
\ge \int_A \int_T u(i, 0,\beta(\gamma(t))) dP^Td\lambda.
\end{eqnarray}
Then, there exists $\delta_1>0$ and a ${\cal I} \boxtimes {\cal T}$ measurable set $D_0 \subseteq A \times T$ such that $\lambda \boxtimes P^T(D_0)>\delta_1$ and for any $(i, t) \in D_0$, $y^p(i,t) \neq 0$ on $D_0$. Thus there exists $\delta>0$, $k \in \{1, \dots, m\}$ and $D \subseteq D_0$ such that $\lambda \boxtimes P^T(D)>\delta$ and the $k$-th component of $y^p(i,t) >\delta$ for any $(i,t) \in D$. Let $e_k$ be the vector in $\mathbb{R}^m_+$ which is one at the $k$-th component and zero otherwise.
Without loss of generality, we can assume $\lambda(D_{t}) < 1$ for any $t \in T$. For $w \in [0,\delta]$, let
\begin{equation*}
\tilde{y}^p(i,t,w)=
\begin{cases}
y^p(i,t)-we_k & \text{if } (i,t) \in D\\
y^p(i,t) & \text{if } (i,t) \notin D,
\end{cases}
\end{equation*}
and $\Phi(i)=\{w: w \in [0,\delta] \text{ and } U_i(\tilde{y}^p_{(i,w)})\in [U_i(y^p_i)-{\frac{1} {4}} d_i,U_i(y^p_i)]\}.$
It is clear that for any fixed $i \in I$, $U_i(\tilde{y}^p_{(i,w)})$ is continuous and decreasing in $w \in [0,\delta]$, and for any fixed $w \in [0,\delta]$, $U_i(\tilde{y}^p_{(i,w)})$ is ${\cal I}$-measurable $i \in I$. Then, $\Phi$ is a compact valued measurable correspondence. For $i \in A$, if $w$ is positive and small enough, the continuity $U_i(\tilde{y}^p_{(i,w)})$ implies that $w \in \Phi(i)$. Let $\phi(i)=max_{w \in \Phi(i)} w$. Then $\phi(i) > 0$ for any $i \in A$. Define an allocation
\begin{equation*}
z^p(i,t)=
\begin{cases}
y^p(i,t)-\phi(i)e_k & \text{if } (i,t) \in D\\
y^p(i,t)+\frac{\int_{D_t}\phi(i)d\lambda}{1-\lambda(D_t)}e_k & \text{if } (i,t) \notin D.
\end{cases}
\end{equation*}
It is clear that the allocation $z^p$ is feasible.

Fix any $i \in A$. It is easy to see for any $t \in T$, $z^p(i,t) \geq \tilde{y}(i,t,\phi(i))$ and $U_i(\tilde{y}(i,\phi(i))) \ge U_i(y^p_i)-{\frac{1} {4}} d_i > U_i(y^p_i)- d_i =U_i(x^p_i)$; then $U_i(z^p_i)> U_i(x^p_i)$.

Next, consider $i \notin A$. Then $(i,t) \notin D$ for any $t \in T$, and thus $z^p_i(t)=x^p_i(t)+\frac{\int_{D_t}\phi(i)d\lambda}{1-\lambda(D_t)}e_k$. Let $D^T=\{t: \lambda(D_t)>0\}$; then $\lambda(D^T)>0$. For any $t \in D^T$, $D_t$ is a subset of $A$ with $\lambda(D_t)>0$, which implies that $\int_{D_t}\phi(i)d\lambda > 0$; hence,
$z^p_i(t)>x^p_i(t)$. For any $t \notin D^T$, we have $\lambda(D_t)=0$, which means that $z^p_i(t)=y^p_i(t)$. The monotonicity assumption on the utility function $u_i$ implies that $U_i(z^p_i)>U_i(y^p_i)$. Therefore, $U_i(z^p_i)>U_i(x^p_i)$ for $\lambda$-almost all $i \in I$, which contradicts the ex ante efficiency of $x^p$. Thus, $\lambda (A)=0$. By (\ref{eq-Jensen}), $x^p$ and $y^p$ are utility equivalent.

(2) \ The result in part (1) says that $U_i (y^p_i)=U_i (x^p_i)$ for $\lambda$-almost all $i \in I$.  The the inequality in (\ref{eq-Jensen}) becomes an equality for $\lambda$-almost all $i \in I$. Assumption 2$^\prime$ requires $u_i$ to be strictly concave in the consumption variable $x\in \mathbb{R}^m_+$. Thus, a strict version of Jensen�s inequality implies that for $\lambda\boxtimes P$-almost all $(i, \omega) \in I \times \Omega$, $\rho(i, \omega)$ is the Dirac measure $\delta_{\int_{\mathbb{R}^m_+}x \; d \rho(i, \omega) }$ at the point $\int_{\mathbb{R}^m_+}x \; d \rho(i, \omega)$.

For any $j \in \{1, \dots, m\}$ and $x\in \mathbb{R}^m_+$, let $\psi_j(x) = \sqrt{x_j}$.
Then, we have
\begin{equation} \label{eq-final}
\mathbb{E}\left(\psi_j(x^p_i)|{\cal C}) \right)=\int_{\mathbb{R}^m_+}\psi_j(x)\ d\delta_{\int_{\mathbb{R}^m_+}x \; d \rho(i, \omega)}  =\psi_j\left(\int_{\mathbb{R}^m_+}x \; d \rho(i, \omega) \right)=\psi_j(\mathbb{E}(x^p_i|{\cal C})).
\end{equation}
For simplicity, let $h_j$ be the function such that $h_j(\omega) = \psi_j(\bar x^p_i(\omega))$ for any $\omega \in \Omega$. Then equation (\ref{eq-final}) implies that $\mathbb{E}\left( h_j |{\cal C}) \right) = \sqrt{\mathbb{E}\left( h^2_j |{\cal C})  \right)}$. That is, $\mathbb{E}\left( h^2_j |{\cal C}  \right) = \left[ \mathbb{E}\left( h_j |{\cal C} \right) \right]^2$, which means that $\mathbb{E}\left( \left( h_j - \mathbb{E}\left( h_j |{\cal C} \right) \right)^2 |{\cal C}  \right) = 0$. Hence,
$\mathbb{E}\left( \left( h_j - \mathbb{E}\left( h_j |{\cal C} \right) \right)^2 \right) = 0$, which implies that $h_j = \mathbb{E}\left( h_j |{\cal C} \right)$. It means that the $j$-th component of $x^p_i(\cdot)$ is ${\cal C}$-measurable for any $j \in \{1, \dots, m\}$. Therefore, we can claim that for $\lambda$-almost all $i \in I$, $x^p_i = \mathbb{E}\left( x^p_i  |{\cal C} \right) = y^p_i$, and we are done.  \qed

\bs

\end{document}